\documentclass[11pt]{amsart}
\usepackage{geometry}                
\usepackage{graphicx}
\usepackage{amssymb}
\usepackage{epstopdf}
\title{On the anisotropic hyperdissipative Navier-Stokes equations}
\author{X-J Wang}

\date{}                                           

\begin{document}
\maketitle

\section{Introduction}
In this note, we consider the global well-posedness of a generalized incompressible Navier-Stokes system in 3D whole space
$$
u_t+u\cdot\nabla u+\nabla p= \nu\mathcal{A}_h u,
$$
\begin{equation}\label{main0}
\nabla\cdot u=0, 
\end{equation}
$$
u(x,0)=u_0(x),
$$
where $u=(u_1, u_2, u_3)\in\mathbf{R}^3$ and $ p$ are the fluid velocity field and pressure. We keep the viscosity constant $\nu=1$ throughout the paper. The initial data $u_0(x)$ is assumed to be smooth, rapidly decreasing and divergence free.  Here $\mathcal{A}_h$  is the anisotropic hyperdissipative operator with
$$\mathcal{A}_hu=-\left( \begin{array}{c} M_1^{2\beta} u_1+M_2^{2\beta} u_1+M_3^{2\alpha} u_1 \\
M_1^{2\beta} u\mathfrak{F}_2+M_2^{2\beta} u_2+M_3^{2\alpha} u_2 \\
   M_1^{2\gamma} u_3+M_2^{2\gamma} u_3+M_3^{2\alpha} u_3 
   \end{array} \right),$$
where $\alpha, \beta, \gamma$ are real numbers and $M_2^{2\alpha}, M_i^{2\beta}, M_i^{2\gamma}, i=1,2 $ are the Fourier multipliers with indices $\alpha, \beta, \gamma$ respectively. For instance, given $u$ in Schwartz space,  $\mathfrak{F}[M_3^{2\alpha} u]=|\xi_3|^{2\alpha} \mathfrak{F}[u]$. Here $\mathfrak{F}$ represents the Fourier transform operator.

When $\alpha=\beta=\gamma=1, \mathcal{A}_h u=\Delta u$, we recover the Navier-Stokes equation.  The  regularity of the Navier-Stokes equation is a consequence of the competition between the nonlinear term $u\cdot\nabla u$ and the diffusion terms $\Delta u$. While the global solutions of Navier-Stokes equations exist in 2D space \cite{Lady59},  it is well-known that the 3D case is still totally open. The difficulty has been interpreted as the diffusive effect being not strong enough to control the nonlinear convective terms. In this aspect, there are naturally two ways to approach the problem. One way is to regularize the convective effect, see for instance the Leray-$\alpha$ model by Cheskidov etc. \cite{ChHO05}. Another way is to strengthen the diffusive effect by increase the Fourier multiplier index $\alpha, \beta, \gamma$ as above. When $\alpha=\beta=\gamma={5/4}$, we have the critical hyperdissipative Navier-Stokes equation, which is well-known to be globally well-posed. See \cite{KaN02} for a proof for the subcritical case in which $\alpha=\beta=\gamma>{5/4}$. The proof in \cite{KaN02} actually applies to the critical case also. A different proof is given in the appendix of Wu\cite{Wu11}, where the author shows the global regularity for the generalized magnetohydrodynamic equations. Our method in this work here suggest another proof, which is presented in \cite{Wang01}. Till now it is an open problem to improve the results from critical ($\alpha=\beta=\gamma={5/4}$) to supercritical ($\alpha=\beta=\gamma<{5/4}$) case.  Tao\cite{Tao09} improves this result slightly towards the supercritical region by establishing global regularity for a logarithmically supercritical hyper-dissipative NS equation, namely for $\mathcal{A}_h$ with $\mathcal{F}[\mathcal{A}_hu]={|\xi|^{5/2}\over g^2(|\xi|)}$ and $\int_1^\infty {ds\over sg(s)^4}=\infty$. Further improvements on the cases of isotropic indices are in great need, but few results are available at this moment for case $1\leq \alpha=\beta=\gamma<{5/4}$, see\cite{CaK82, KaN02}.

Starting with the Navier-Stokes equation, a natural question arises: what if we do not raise all the indices from $1$ to $5/4$? To get the global well-posedness, is it enough to increase the indices of some components of the Laplacian while keep the others? In other words, we are interested in the structure of the anisotropic hyperdissipative(fractional Laplacian) operator.

Suppose $\alpha\geq 3/2, \beta=1, \gamma\geq 5/4, $, namely when

\begin{equation}\label{op1}
\mathcal{A}_hu= \left( \begin{array}{c} \partial_{x_1x_1} u_1+\partial_{x_2x_2} u_1-M_3^{2\alpha} u_1 \\
\partial_{x_1x_1} u_2+\partial_{x_2x_2} u_2-M_3^{2\alpha} u_2 \\
   -M_1^{2\gamma} u_3-M_2^{2\gamma} u_3-M_3^{2\alpha} u_3 
   \end{array} \right),
\end{equation}

We have our main result as following.

\subsection{Theorem 1.1.}
With $\mathcal{A}_hu$ given  in (\ref{op1}), the anisotropic generalized Navier-Stokes system (\ref{main0}) has a global smooth solution.

\subsection{Remarks} We make a few remarks.

1. This is certainly an artificial model, just like the isotropic hyperdissipative model with $\mathcal{A}_hu=-(-\Delta)^{5\over4}u$. However, they provide different perspectives to approximating the Navier-Stokes equations, given the difficulties in the study of this problem \cite{Taoblog}. To certain extent, it is actually a common practice to inject artificial diffusion into the system in both the analytical and the numerical study, see for instance \cite{SaT76, Jim94}.

2. The analytical study of the hyperdissipative operator itself is of certain interests. Typically, using techniques in \cite{Tao09}, we can make a logarithmical improvement on our result. But we refrain us from doing that for sake of clarity. 

3. By adding partial hyperdissipation, the system loses the well-known scale invariance $u^\lambda={1\over\lambda}u({x\over\lambda},{t\over{\lambda^2}})$. The breakdown of scale invariance seems to be unavoidable in constructing a quantity that bounds the energy, see Wang\cite{Wang01}. Indeed, most regularization schemes used in PDE analysis seem to cause breakdown of the scale invariance.

4. Upon the submission of this work, the authors notice that Zhang\cite{Zhang10} has worked on a similar problem using Littlewood-Paley theory. But unfortunately, there are a couple of flawed calculations and their conclusion seems to be incorrect.
The system considered in \cite{Zhang10} does obey a scale invariance, but under that scale transformation, $\sup\limits_{0\leq t\leq T}||\nabla u||_{L^2}$ is a supercritical quantity. The effort to use such a quantity to stop blowup is doomed to fail, see the analysis in\cite{Wang01}.

5. A weaker version the theorem would be the case where 
$$\mathcal{A}_hu= \Delta u-\left( \begin{array}{c} M_3^{2\alpha} u_1 \\
M_3^{2\alpha} u_2 \\
M_1^{2\gamma} u_3+M_2^{2\gamma} u_3+M_3^{2\alpha} u_3 
\end{array} \right).$$ It can be proved using the techniques here with minor modifications. Comparing to (\ref{op1}), it is easier to implement the numerical regularization. 

6. Merely for clarity purpose, we will assume $\gamma\leq \alpha$ in the proof. There would little difference in the proof otherwise.
\subsection{Notations}

We consider 3D space in this note.  We use $\nabla_h u$ for $(\partial_{x_1}u, \partial_{x_2} u)$. The Sobolev $L^2$ norm on the whole space is $||\cdot||_{L^2(\mathbf{R}^3)}$, or $||\cdot ||_{L^2}$ for short. Similar for $L^p$ norm $||\cdot||_{L^p}$ and $W^{s,p}$  norm $||\cdot||_{W^{s,p}}$.We use $||\cdot||_{L_h^2}$ and  $||\cdot||_{L_v^\infty}$ for the $L^2$ norm in $x_1x_2$ plane and $L^\infty$ norm in $x_3$ direction, respectively. We use $X\lesssim Y$ to denote the statement that $X\lesssim CY$ for some absolute constant $C>0$.
The Fourier transform is defined by $$\mathfrak{F}[u](\xi)={1\over (2\pi)^{3\over 2}}\int_{\mathbf{R}^3}e^{-i\xi x}u(x)dx,$$ 
which is often denoted by $\hat{u}$ later on. We use $<\cdot, \cdot>$ to denote the $L^2$ pair of two functions in either the physical space or frequency space.
We will use Parseval's theorem to frequently transform an integral between its physical space and frequency space. For instance we have $||u||_{L^2_x}=||\hat{u}||_{L^2_{\xi}}$ and $$<\partial_{x_1x_1} u_1+\partial_{x_2x_2} u_1-M_3^{2\alpha} u_1, u_1>=\int (\xi_1^2+\xi_2^2+\xi_3^{2\alpha})\hat{u}_1^2d\xi$$
$$=||\sqrt{\xi_1^2+\xi_2^2}\hat{u}_1||_{L^2}^2+|||\xi_3|^\alpha u_1||^2_{L^2}=||\nabla_h u_1||^2+||M_3^\alpha u_1||^2.$$

\section{Lemmas} We prepare a list of lemmas.

{\bf Lemma 2.1.}  For any function $u(x): \mathbf{R}^3\to \mathbf{R}$, we have 
\begin{equation}
||u||_{L_h^4L_v^2}\leq ||u||_{L_v^2L_h^4}\lesssim||u||_{L^2}^{1\over 2}||\nabla_h u||_{L^2}^{1\over 2},
\end{equation} provided the right hand side exists.
Here the subscript $h$ represents the horizontal $x_1x_2$ plane and $v$ represents the vertical $x_3$ direction. For the proof, see Paicu\cite{Pa05}.

We need an embedding inequality of logarithmical type due to Brezis \cite{BreGa80}.

{\bf Lemma 2.2.}  For smooth function $\phi(x): \mathbf{R}\to \mathbf{R}$ in $H^1(\mathbf{R})$, we have 
$$
||\phi(x)||_{L^\infty(\mathbf{R})}\lesssim \sqrt{\ln(1+N)}||\phi||_{H^{1\over 2}(\mathbf{R})}+{1\over\sqrt{1+N}}||\phi||_{H^1(\mathbf{R})}
$$

Here $N>0$ is a number to be determined later in the proof of our theorem. 

{\it Proof.} By Fourier transform $|\phi(x)|=|\int e^{i\xi x}\hat{\phi}(\xi)d\xi|\leq\int |\hat{\phi}(\xi)|d\xi$ and splitting the right hand side into
$\int_{|\xi|
\leq N}|\hat{\phi}(\xi)|(1+|\xi|)^{1\over 2}{1\over (1+|\xi|)^{1\over 2}}d\xi+\int_{|\xi|> N}|\hat{\phi}(\xi)|(1+|\xi|){1\over (1+|\xi|)}d\xi$. By Cauchy-Schwarz inequality, a straightforward computation leads to the conclusion.

{\bf Corollary 2.3.} Let $\alpha\geq {3\over 2}$. We have for smooth function $\phi(x)$
\begin{equation}\label{Bre}
||\partial_x\phi(x)||_{L^\infty(\mathbf{R})}\lesssim \sqrt{\ln(1+N)}(||\phi||_{L^ 2}+||M^\alpha \phi||_{L^2})+{1\over\sqrt{1+N}}(||\phi||_{L^ 2}+||M^{\alpha+1}\phi||_{L^2}).
\end{equation}

Here $M^\alpha$ is of course the 1D Fourier multiplier with respect to $x$.

{\it Proof.}
Apply {\bf Lemma 2.2.} to $\partial_x\phi(x)$. Note that $|\xi|(1+|\xi|)^{1\over 2}\lesssim (1+|\xi|^{\alpha}),$ and $|\xi|(1+|\xi|)\lesssim (1+|\xi|^{\alpha+1}), \alpha\geq {3/2}$.

%
%

We also need an inequality due to Ladyzhenskaya \cite{Lady59}.

{\bf Lemma 2.4.} In 3D space, we have 
\begin{equation}\label{lady}
||\phi(x)||_{L^4}^2\lesssim ||\phi||_{L^2}^{1\over 2}||\nabla \phi||_{L^2}^{3\over 2}
\end{equation}

\section{Proof of Theorem 1.1.}

Consider
\begin{eqnarray}\label{main2}
u_t+u\cdot\nabla u+\nabla p&=&  \left( \begin{array}{c} \partial_{x_1x_1} u_1+\partial_{x_2x_2} u_1-M_3^{2\alpha} u_1 \\
\partial_{x_1x_1} u_2+\partial_{x_2x_2} u_2-M_3^{2\alpha} u_2 \\
   -M_1^{2\gamma} u_3-M_2^{2\gamma} u_3-M_3^{2\alpha} u_3 
   \end{array} \right),
\end{eqnarray}  
$$\nabla\cdot u=0,$$
$$u(x,0)=u_0(x).$$

The local existence proof is a simple application of energy estimate, same as the proof for Navier-Stokes equation. More specifically, for any data $u_0\in H^s, s>{5\over 2}$, there exists $T>0$, depending on $||u_0||_H^s$, so that the solution to (\ref{main2}) exists in $C([0,T], H^s)\cap AC([0,T], H^{s-1})$, see\cite{Kato72}. 

To achieve the global existence, we show that on any time interval $[0, T)$ where the solution exists, $||u||_{H^s}$ stays bounded as $t$ approaches $T$. No blowup of $||u||_{H^s}$ occurring for any such $T$ implies global existence of smooth solutions to (\ref{main2}). We first get $L^2$ and $H^1$ estimates, then we bound the $H^s$ norm of the solution for any $s>5/2$ and initial data $u_0\in H^s$. Since we assume smooth initial data $u_0$, the solution is smooth as well.

\subsection{$L^2$ Energy estimate} 

Multiplying (\ref{main2}) by $u$, using incompressible condition to get rid of $\nabla p$, we have 

\begin{equation}\label{foru}
{1\over 2}{d\over dt}||u||_{L^2(\mathbf{R}^3)}^2 =-m(t),
\end{equation}
with $ m(t)= m_1(t)+m_2(t)+m_3(t)$, and

$$m_1(t)=||\nabla_h u_1||^2+||\nabla_h u_2||^2,$$
$$m_2(t)= ||M_1^\gamma u_3||^2+||M_2^\gamma u_3||^2,$$
$$m_3(t)=||M_3^\alpha u_1||^2 +||M_3^\alpha u_2||^2 +||M_3^\alpha u_3||^2.$$


So we have $$\sup\limits_{0\leq t\leq T}||u||_{L^2}\leq C,$$
$$\int_0^T m(t) dt\leq C.$$

In particular $$\int_0^T ||u||_{L^2}+ m(t) dt\leq C.$$

\subsection{$H^1$ Energy estimate}

Multiplying (\ref{main2}) by $\Delta u$, we have

\begin{equation}\label{grad1}
{1\over 2}{d\over dt}||\nabla u||^2 + \tilde{m}(t)=\int u\cdot\nabla u\cdot \Delta udx
\end{equation}

with $ \tilde{m}(t)= \tilde{m}_1(t)+\tilde{m}_2(t)+\tilde{m}_3(t)$, and

$$\tilde{m}_1(t)=||\nabla_h \nabla u_1||^2+||\nabla_h\nabla u_2||^2,$$
$$\tilde{m}_2(t)= ||M_1^\gamma \nabla u_3||^2+||M_2^\gamma\nabla  u_3||^2,$$
$$\tilde{m}_3(t)=||M_3^\alpha \nabla u_1||^2 +||M_3^\alpha\nabla  u_2||^2 +||M_3^\alpha\nabla  u_3||^2.$$

The right hand side of (\ref{grad1}) is equal to 

\begin{eqnarray*}
\sum_{i,j,k=1}^3\int u_i\partial_i u_j \partial_{kk}u_jdx=-\sum_{i,j,k=1}^3\int \partial_ku_i\partial_i u_j \partial_{k}u_jdx-\sum_{i,j,k=1}^3\int u_i\partial_i \partial_ku_j \partial_{k}u_jdx.
\end{eqnarray*}

The second term on the right disappears because of the incompressible condition. 
So we have

\begin{eqnarray}
{1\over 2}{d\over dt}||\nabla u||^2 + \tilde{m}(t)=-\sum_{i,j,k=1}^3\int \partial_ku_i\partial_i u_j \partial_{k}u_jdx.
\end{eqnarray}

Now we expect $\tilde{m}(t)$ to control the convective terms on the right. We know in case of $\alpha=\gamma=1$, hence $\mathcal{A}_hu=\Delta u$, the dissipation is weak and $\tilde{m}(t)$ is not strong enough to dominate the nonlinear convective terms. On the other hand, if $\mathcal{A}_hu=-(-\Delta)^{5\over 4} u$, it is a standard result \cite{KaN02} that the dissipation is strong enough to control the nonlinear effects and achieve the global regularity. Now in our the case, the difficulty is clearly coming from the fact that only the right and bottom marginal components in our hyper-dissipative operators have strong dissipation. We have to use the marginal strong hyper-dissipation to control the full nonlinear terms.

We regroup the summations 
\begin{eqnarray}
\sum_{i,j,k=1}^3\int \partial_ku_i\partial_i u_j \partial_{k}u_jdx=A+B
\end{eqnarray}
where 

$$
A=\sum_{i,j,k=1}^2\int \partial_ku_i\partial_i u_j \partial_{k}u_jdx
$$

while $B$ contains all the terms $ \int \partial_ku_i\partial_i u_j \partial_{k}u_jdx$ with either $i=3, j=3$ or $k=3$.

By Lemma2.2 in Kukavica etc. \cite{KuZ07},
\begin{eqnarray*}
A=\sum_{i,j,k=1}^2\int \partial_ku_i\partial_i u_j \partial_{k}u_jdx=-\sum_{i,j=1}^2\int \partial_iu_j\partial_i u_j \partial_{3}u_3dx
\end{eqnarray*}

\begin{eqnarray*}
+\int\partial_1u_1\partial_2u_2\partial_3u_3dx-\int\partial_1u_2\partial_2u_1\partial_3u_3dx
\end{eqnarray*}

A key observation is that each term in $A$ contains $\partial_3 u_3$ while each term in $B$ contains either $\partial_3 u_i$ or $\partial_i u_3\partial_j u_3$.  All the terms in $A$ and $B$ are the `marginal components' in the hyper-dissipative term $\mathcal{A}_hu$. They belong to two groups.   Terms containing $\partial_3 u_i$ belong to group 1, we bound them by splitting the space into horizontal and vertical coordinates and using Lemma1. We put terms containing $\partial_i u_3\partial_j u_3, i\neq 3, j\neq 3$ into group 2, which can be controlled by Gagliardo-Nirenberg inequality. Now the strategy for the proof of the theorem is evident. We use the marginal components to bound the full set of components. We pick one term in each group to show the idea. Proofs of the rest are the same.

\subsubsection{group1, terms like $\int \partial_iu_j\partial_ku_l\partial_3u_mdx$}

$$
\int \partial_iu_j\partial_ku_l\partial_3u_mdx\lesssim\int |\nabla u|^2|\partial_3 u|dx\lesssim\int\int \{\int |\nabla u|^2|\partial_3 u|dx_3\}dx_1dx_2
$$

$$
\lesssim \parallel||\nabla u||^2_{L_v^2}||\partial_3 u||_{L_v^{\infty}}\parallel_{{L_h^1}}
$$

Thanks to (\ref{Bre}) in {\bf Corollary 2.3.}, we have $\int \partial_iu_j\partial_ku_l\partial_3u_mdx$

$$
\lesssim \parallel||\nabla u||^2_{L_v^2}\Big(\sqrt{\ln(1+N)}(||u||_{L^2_v}+||M_3^\alpha u||_{L_v^{2}})+{1\over \sqrt{1+N}}(||u||_{L^2}+||M_3^\alpha \partial_3 u||_{L_v^{2}})\Big)\parallel_{L_h^1}
$$

$$
\lesssim \parallel||\nabla u||^2_{L_v^2}\parallel_{L_h^2}\parallel \Big(\sqrt{\ln(1+N)}(||u||_{L^2_v}+||M_3^\alpha u||_{L_v^{2}})+{1\over \sqrt{1+N}}(||u||_{L^2}+||M_3^\alpha \partial_3 u||_{L_v^{2}})\Big)\parallel_{L_h^2}
$$

\begin{equation*}
\lesssim ||\nabla u||^2_{L_h^4L_v^2}\Big(\sqrt{\ln(1+N)}(||u||_{L^2}+||M_3^\alpha u||_{L^{2}})+{1\over \sqrt{1+N}}(||u||_{L^2}+||M_3^\alpha \partial_3 u||_{L^{2}})\Big)
\end{equation*}

By {\bf Lemma 2.1.}, we have $$\int \partial_iu_j\partial_ku_l\partial_3u_mdx$$
\begin{equation*}
\lesssim ||\nabla u||_{L^2}||\nabla_h\nabla u||_{L^2}\Big(\sqrt{\ln(1+N)}(||u||_{L^2}+||M_3^\alpha u||_{L^{2}})+{1\over \sqrt{1+N}}(||u||_{L^2}+||M_3^\alpha \partial_3 u||_{L^{2}})\Big)
\end{equation*}

\begin{equation*}
\leq  \delta||\nabla_h\nabla u||^2_{L^2}+C||\nabla u||^2_{L^2} \Big(\ln(1+N)(||u||_{L^2}^2+||M_3^\alpha u||^2_{L^{2}})+{1\over 1+N}(||u||_{L^2}^2+||M_3^\alpha \partial_3 u||_{L^{2}}^2)\Big)
\end{equation*}
where $\delta$ is to be determined soon.

We now work on $||\nabla_h\nabla u||_{L^2}^2$. Recall the definition of $\tilde{m}(t)$ in (\ref{grad1}). we have 
$$||\nabla_h\nabla u||_{L^2}^2=\tilde{m}_1(t)+||\nabla_h\nabla u_3||_{L^2}^2,$$
and since $\gamma\geq {5/4}$, we have
$$||\nabla_h\nabla u_3||_{L^2}^2=\int (\xi_1^2+\xi_2^2)|\xi|^2\hat{u}_3^2d\xi\lesssim \int (1+\xi_1^{2\gamma}+\xi_2^{2\gamma})|\xi|^2\hat{u}_3^2d\xi\lesssim \tilde{m}_2(t)+||\nabla u_3||_{L^2}^2.$$

So $$||\nabla_h\nabla u||_{L^2}^2\leq C(\tilde{m}_1(t)+\tilde{m}_2(t)+||\nabla u_3||_{L^2}^2).$$ 

Picking $\delta$ so that $\delta C\leq 1/54$, we have 

$$\int \partial_iu_j\partial_ku_l\partial_3u_mdx$$

\begin{equation*}
\leq  {1\over 54}[\tilde{m}_1(t)+\tilde{m}_2(t)]+C||\nabla u||^2_{L^2} \Big(\ln(e+N)(||u||_{L^2}^2+||M_3^\alpha u||^2_{L^{2}})+{1\over 1+N}(||u||_{L^2}^2+||M_3^\alpha \partial_3 u||_{L^{2}}^2)+1\Big)
\end{equation*}

Now pick $N= 54C||\nabla u||_{L^2}^2$, we have $\int \partial_iu_j\partial_ku_l\partial_3u_mdx$
\begin{equation*}
\leq {1\over 54}[\tilde{m}_1(t)+\tilde{m}_2(t)]+{1\over 54}||M_3^\alpha \partial_3 u||^2_{L^2}+ C||\nabla u||^2_{L^2} \ln(e+||\nabla u||^2_{L^2} )\Big(||M_3^\alpha u||^2_{L^{2}}+||u||_{L^2}^2+1\Big)
\end{equation*}

\begin{equation*}
\leq {1\over 54}\tilde{m}(t)+ C||\nabla u||^2_{L^2} \ln(e+||\nabla u||^2_{L^2} )\Big(||M_3^\alpha u||^2_{L^{2}}+||u||_{L^2}^2+1\Big)
\end{equation*}

\begin{equation*}
\leq {1\over 54}\tilde{m}(t)+ C||\nabla u||^2_{L^2} \ln(e+||\nabla u||^2_{L^2} )\Big(m(t)+||u||_{L^2}^2+1\Big)
\end{equation*}

\begin{equation}\label{group1}
\leq {1\over 54}\tilde{m}(t)+ C(e+||\nabla u||^2_{L^2}) \ln(e+||\nabla u||^2_{L^2} )\Big(m(t)+||u||_{L^2}^2+1\Big)
\end{equation}

This finishes the estimates for terms in group 1. 

\subsubsection{group2} Terms like $\int \partial_iu_j\partial_ku_3\partial_lu_3dx, i\neq 3, j\neq 3$.  can be bounded by

$$\int |\nabla u||\nabla u_3||\nabla u_3|dx\lesssim ||\nabla u||_{L^2}||u_3||_{W^{1,{12\over 5}}}||u_3||_{W^{1,12}}\lesssim ||\nabla u||_{L^2}||u_3||_{H^\gamma}||u_3||_{H^{\gamma+1}}
$$

$$\leq \delta' ||u_3||_{H^{\gamma+1}}^2+C||\nabla u||^2_{L^2}||u_3||^2_{H^\gamma},$$
where $\delta'$ is a small positive number to be determined.

Noting that $\gamma\leq\alpha$, we have 
$$||u_3||^2_{H^\gamma}\lesssim \int(1+\xi_1^{2\gamma}+\xi_2^{2\gamma}+\xi_3^{2\alpha})\hat{u}^2_3d\xi\lesssim ||u||^2_{L^2}+m_3(t)+m_2(t)\lesssim ||u||^2_{L^2}+m(t).$$

Moreover, 
$$||u_3||^2_{H^{\gamma+1}}\lesssim \int (1+[(1+\xi_1^{2\gamma}+\xi_2^{2\gamma}+\xi_3^{2\alpha})|\xi|^2])\hat{u}^2_3d\xi$$
$$\lesssim ||u||^2_{L^2}+||\nabla u||^2_{L^2}+\tilde{m}_3(t)+\tilde{m}_2(t)\lesssim ||u||^2_{L^2}+||\nabla u||^2_{L^2}+\tilde{m}(t).$$

We can pick $\delta'$ so small that $$\int |\nabla u||\nabla u_3||\nabla u_3|dx$$
$$\leq {1\over 54} \tilde{m}(t)+C||\nabla u||^2_{L^2}\{1+||u||^2+m(t)\}+C||u||^2_{L^2}$$

\begin{equation}\label{group2}
\leq {1\over 54} \tilde{m}(t)+C(e+||\nabla u||^2_{L^2})\{1+||u||^2+m(t)\}\end{equation}

\subsubsection{Completion of $H^1$ estimate}
Combining the estimate (\ref{group1})-(\ref{group2}), noting there are totally 27 terms, we have
\begin{equation*}
{d\over dt}||\nabla u||^2 +\tilde{m}(t)
\lesssim (e+||\nabla u||^2_{L^2}) \ln(e+||\nabla u||^2_{L^2} )\{1+||u||^2+m(t)\}
\end{equation*}

Since $\{1+||u||^2+m(t)\}$ is integrable on $[0,T]$, this gives a double exponential bound for $(e+||\nabla u||^2_{L^2})$, hence for $||\nabla u||^2_{L^2}$. More specifically, we have 
\begin{equation*}
\int_0^T\tilde{m}(s)ds\leq C,
\end{equation*}
\begin{equation}\label{grad2}
\sup\limits_{0\leq t\leq T}||\nabla u||\leq C.
\end{equation}

\subsection{High order energy estimate}

Let $E_s$ denote the energy, $$E_s=||\nabla^s u||_{L^2}^2+||u||^2_{L^2},$$ which is an equivalent norm to $||u||_{H^s}$. We establish a differential inequality to bound $E_s$.

Multiplying (\ref{main2}) by $(-\Delta)^su$, integrating on space and combining with (\ref{foru}), we have 
\begin{equation}\label{high1}
{1\over 2} {d\over dt} E_s+ m(t) + \hat{m}(t)=-\int u\cdot\nabla u\cdot (-\Delta)^s udx=-\int \nabla^s(u\cdot\nabla u)\cdot \nabla^s udx 
\end{equation}


with $ \hat{m}(t)= \hat{m}_1(t)+\hat{m}_2(t)+\hat{m}_3(t)$, and

$$\hat{m}_1(t)=||\nabla_h \nabla^s u_1||^2+||\nabla_h\nabla^s u_2||^2,$$
$$\hat{m}_2(t)= ||M_1^\gamma \nabla^s u_3||^2+||M_2^\gamma \nabla^s u_3||^2,$$
$$\hat{m}_3(t)=||M_3^\alpha \nabla^s u_1||^2 +||M_3^\alpha \nabla^s u_2||^2 +||M_3^\alpha \nabla^s u_3||^2.$$

We rewrite the right hand side of (\ref{high1}) as 

\begin{equation}\label{high2}
\sum_{0\leq \sigma \leq s}C_{\sigma}\int\nabla^\sigma u \nabla^{s-\sigma+1}u\nabla^sudx,
\end{equation}
where $C_{\sigma}$ are the combinatoric coefficients, which are of no importance in our estimate. Because of the incompressibility condition, the term with $\sigma=0$ disappears. The term with $\sigma=s$ equals to 
$$C_{s}\int\nabla u \nabla^{s}u\nabla^s udx \lesssim ||\nabla u||_{L^2}||\nabla^s u||_{L^4}^2\lesssim ||\nabla u||_{L^2}||\nabla^s u||_{L^2}^{1/2}||\nabla^{s+1} u||_{L^2}^{3/2} $$ where we have used Ladyzhenskaya inequality (\ref{lady}). 

As a matter fact, we claim that all the rest terms in (\ref{high2}) can be bound by $\int\nabla u \nabla^{s}u\nabla^s udx$. Indeed, H$\ddot{o}$lder inequality implies 
$$\int\nabla^\sigma u \nabla^{s-\sigma+1}u\nabla^sudx\lesssim ||\nabla^\sigma u||_{L^p} ||\nabla^{s-\sigma+1}u||_{L^q}||\nabla^su||_{L^4}$$ for any integer $0<\sigma<s,$ and for any $p>1,q>1, 1/p+1/q=4/3$.
Applying Gagliardo-Nirenberg inequality to both $||\nabla^\sigma u||_{L^p}$ and $ ||\nabla^{s-\sigma+1}u||_{L^q}$, we have

$$||\nabla^\sigma u||_{L^p}=||\nabla^{\sigma-1}\nabla u||_{L^p}\lesssim ||\nabla^{s-1} \nabla u||_{L^4}^\theta ||\nabla u||_{L^2}^{1-\theta}$$

$$||\nabla^{s-\sigma+1}u||_{L^q}\lesssim||\nabla^{s-1} \nabla u||_{L^4}^{1-\theta} ||\nabla u||_{L^2}^\theta.$$

Combining the above estimates, we  have 
$$\sum_{0\leq \sigma \leq s}C_{\sigma}\int\nabla^\sigma u \nabla^{s-\sigma+1}u\nabla^sudx\lesssim ||\nabla u||_{L^2}||\nabla^s u||^2_{L^4},
$$

Using {\bf Lemma 2.4.} and Young's inequality, we can bound  the right hand side of (\ref{high1}) by
$$\sum_{0\leq \sigma \leq s}C_{\sigma}\int\nabla^\sigma u \nabla^{s-\sigma+1}u\nabla^sudx\lesssim ||\nabla u||_{L^2}||\nabla^s u||_{L^2}^{1/2}||\nabla^{s+1} u||_{L^2}^{3/2},
$$

$$\leq \sigma'' ||\nabla^{s+1}u||^2_{L^2}+C||\nabla u||^4_{L^2}||\nabla^s u||_{L^2}^2$$

We claim that $ ||\nabla^{s+1}u||_{L^2}^2\lesssim \hat{m}(t)+ ||\nabla^su||_{L^2}^2$. Once this is proved, we can choose $\sigma''$ so that 
$$\sum_{0\leq \sigma \leq s}C_{\sigma}\int\nabla^\sigma u \nabla^{s-\sigma+1}u\nabla^sudx\leq {1\over 2} \hat{m}(t)+C(1+||\nabla u||^4_{L^2})||\nabla^s u||_{L^2}^2$$

Hence by (\ref{high1}) we have
$$
{d\over dt} E_s+ 2m(t)+ \hat{m}(t)\leq C(||\nabla u||^4_{L^2}+1)||\nabla^s u||_{L^2}^2.
$$
In particular we have
$$
{d\over dt} E_s\leq C(||\nabla u||^4_{L^2}+1)E_s.
$$
By Gronwall's inequality and the fact $\sup\limits_{0\leq t\leq T}||\nabla u ||_{L^2}\leq C$, the high order estimate is bounded. 
 
We are left with the proof of $ ||\nabla^{s+1}u||_{L^2}^2\lesssim \hat{m}(t)+ ||\nabla^su||_{L^2}^2$.  Actually 

$$||\nabla^{s+1}u||_{L^2}^2=\int |\xi|^{2s+2}|\hat{u}|^2d\xi=\int |\xi|^2|\xi|^{2s}|\hat{u}|^2d\xi$$
$$=\int (|\xi_1|^2+|\xi_2|^2+|\xi_3|^2)|\xi|^{2s}(|\hat{u}_1|^2+|\hat{u}_2|^2+|\hat{u}_3|^2)d\xi$$
$$=\int (|\xi_1|^2+|\xi_2|^2)|\xi|^{2s}(|\hat{u}_1|^2+|\hat{u}_2|^2)d\xi+\int (|\xi_1|^2+|\xi_2|^2)|\xi|^{2s}|\hat{u}_3|^2d\xi+\int |\xi_3|^2|\xi|^{2s}|\hat{u}|^2d\xi$$

The first integral is exactly $\hat{m}_1(t)$. Recall $\gamma\geq 5/4$, the second integral gives $$\int (|\xi_1|^2+|\xi_2|^2)|\xi|^{2s}|\hat{u}_3|^2d\xi$$
$$\lesssim \int (1+|\xi_1|^{2\gamma}+|\xi_2|^{2\gamma})|\xi|^{2s}|\hat{u}_3|^2d\xi$$
$$\lesssim ||\nabla u_3||^2+\hat{m}_2(t)$$

The third integral 
$$\int |\xi_3|^2|\xi|^{2s}|\hat{u}|^2d\xi\lesssim \int (1+|\xi_3|^{2\alpha})|\xi|^{2s}|\hat{u}|^2d\xi\lesssim ||\nabla u||^2+\hat{m}_3(t).$$

Thus we complete the proof of $||\nabla^{s+1}u||_{L^2}^2\lesssim \hat{m}(t)+ ||\nabla^su||_{L^2}^2$.

\end{document}